\documentclass[11pt]{article}
\usepackage{amssymb}
\usepackage{color}

\normalbaselines

\newtheorem{Theorem}{Theorem}[section]

\newtheorem{Proposition}[Theorem]{Proposition}

\newtheorem{Lemma}[Theorem]{Lemma}

\newtheorem{Corollary}[Theorem]{Corollary}

\newtheorem{Remark}[Theorem]{Remark}

\newcommand{\RR}{{{\rm I} \kern -.15em {\rm R} }}

\newcommand{\C}{{{\rm l} \kern -.42em {\rm C} }}

\newcommand{\nat}{{{\rm I} \kern -.15em {\rm N} }}

\newcommand{\be}{\begin{equation}}
\newcommand{\ee}{\end{equation}}
\newcommand{\beq}{\begin{eqnarray}}
\newcommand{\eeq}{\end{eqnarray}}
\newcommand{\beqs}{\begin{eqnarray*}}
\newcommand{\eeqs}{\end{eqnarray*}}
\newcommand{\bt}{\begin{Theorem}}
\newcommand{\et}{\end{Theorem}}
\newcommand{\br}{\begin{Remark}}
\newcommand{\er}{\end{Remark}}
\newcommand{\bc}{\begin{Corollary}}
\newcommand{\ec}{\end{Corollary}}
\newcommand{\bl}{\begin{Lemma}}
\newcommand{\el}{\end{Lemma}}
\newcommand{\bd}{\begin{definition}}
\newcommand{\ed}{\end{definition}}

\title{Stabilization of second--order evolution\\
 equations with
time delay}
\author{{\sc Serge Nicaise}
\\Institut des Sciences et Techniques de Valenciennes
\\Laboratoire de Math\'ematiques et leurs Applications
\\Universit\'e de Valenciennes et du Hainaut Cambr\'esis
\\
59313 Valenciennes Cedex 9, France\\\\
{\sc Cristina Pignotti}
\\Dipartimento di Ingegneria e Scienze dell'Informazione e Matematica\\
 Universit\`{a} di L'Aquila\\
Via Vetoio, Loc. Coppito, 67010 L'Aquila Italy}

\date{}

\begin{document}

\textwidth=160 mm

\textheight=225mm

\parindent=8mm

\frenchspacing

\maketitle

\begin{abstract}
We consider second--order evolution equations  in an abstract setting with damping and time delay
 and  give sufficient conditions ensuring exponential stability. Our abstract framework is then applied to the wave equation,  the elasticity system and the Petrovsky system.
\end{abstract}

\vspace{5 mm}

\def\qed{\hbox{\hskip 6pt\vrule width6pt
height7pt
depth1pt  \hskip1pt}\bigskip}



\section{Introduction}
\label{pbform}\hspace{5mm}

\setcounter{equation}{0}

Let $H$ be a real Hilbert space  with norm and inner product denoted respectively by
 $\Vert \cdot\Vert_{H}$ and $\langle \cdot ,\cdot\rangle_{H}$
and let $A:{\mathcal D}(A)\rightarrow H$
be a positive self--adjoint  operator with a compact inverse in $H.$ Denote by $V:={\mathcal D}(A^{\frac 1 2})$ the domain of
$A^{\frac 1 2}.$ Moreover, for $i=1,2,$ let $U_i$ be real Hilbert spaces
 with norm and inner product denoted respectively by
 $\Vert \cdot\Vert_{U_i}$ and $\langle \cdot ,\cdot\rangle_{U_i}$
 and let $B_i: U_i\rightarrow V'$ be linear operators.
In this setting we consider the problem
 \begin{eqnarray}
& &u_{tt}(t) +A u (t)+B_1B_1^*u_t(t) +B_2B_2^*u_t(t-\tau) =0\quad t>0,\label{1.1}\\
& &u(0)=u_0\quad \mbox{\rm and}\quad u_t(0)=u_1,\quad\label{1.2}\\
& &B_2^*u_t(t)=f^0(t) \quad t\in (-\tau ,0), \label{1.3}
\end{eqnarray}
where the constant  $\tau >0$ is the time delay.
We assume that the delay feedback operator $B_2$ is bounded, that is $B_2\in \mathcal{L}(U_2,H)$, while the standard one  
$B_1\in 
\mathcal{L}(U_1, V^\prime)$
may be   unbounded.

Time delays
are often present
in applications and practical problems and it is by now well--known that even an arbitrarily small delay in the feedback may destabilize a system which is uniformly exponentially stable in absence of delay.
For some examples in this sense we refer to \cite{Datko, DLP, NPSicon06, XYL}.

We are interested in giving     stability results
for the above problem under a suitable assumption on the ``size '' of the feedback operator $B_2,$  when the feedback $B_1$ is a stabilizing one.
More precisely, we will show that for a system which is exponentially stable in absence of time delay, i.e. for
$B_2=0,$ the exponential stability is preserved if $\Vert B_2^*\Vert$
is sufficiently small.

In this sense this paper extends and generalizes the result of \cite{pignotti}
for wave equation with local damping and time delay. On the other hand it
completes the analysis of \cite{NPSicon06, NV2010}.
Indeed here we do not assume
$$\exists\ \alpha <1 \ \mbox{such that}\ \Vert B_2^*u\Vert_{U_2}\le \alpha \Vert B_1^*u\Vert_{U_1},\ \forall\ u\in V;$$
as in \cite{NV2010} (cfr. assumption (1.8) of \cite{NPSicon06} for the wave equation).

Assuming that an observability inequality holds for the system (\ref{1.1}), (\ref{1.2}) when $B_2=0,$ through the definition of a suitable energy (see (\ref{energy}))
and the use of a perturbation argument as in \cite{pignotti}, we obtain sufficient conditions ensuring exponential stability.
Our abstract framework is then applied to some concrete examples, namely the wave equation, the elasticity system and the Petrovsky system.

The paper is organized as follows. In section \ref{well}    a well--posedness result of the abstract system is proved. In section \ref{st} we obtain exponential stability results for the abstract system under suitable conditions. Finally, in sections \ref{WW}, \ref{Elastic} and
\ref{PP}
we apply our abstract results to the wave equation with local and boundary dampings,  the elasticity system and the Petrovsky system respectively. Other examples (like wave or beam equations on networks)
could be given, we skip them for shortness.

\section{Well-posedness \label{well}}

\hspace{5mm}

\setcounter{equation}{0}

In this section we will give well--posedness results for problem
(\ref{1.1})--(\ref{1.3})
using semigroup theory.

As in \cite{NPSicon06} we introduce the function
\begin{equation}\label{zeta}
z(\rho ,t):=B_2^* u_t(t-\tau\rho ),\quad \rho\in (0,1),\ t>0.
\end{equation}

Then, problem (\ref{1.1})--(\ref{1.3}) can be rewritten as

 \begin{eqnarray}
& &u_{tt}(t) +A u (t)+B_1B_1^*u_t(t) +B_2z(1,t) =0\quad t>0,\label{N.1}\\
& & \tau z_t(\rho ,t)+z_\rho (\rho ,t) =0,\quad \rho\in (0,1),\ t>0, \label{Z.1}\\
& &u(0)=u_0\quad \mbox{\rm and}\quad u_t(0)=u_1,\quad\label{N.2}\\
& &z(\rho , 0)=f^0(-\tau\rho) \quad \rho \in (0,1), \label{N.3}\\
& &z(0, t)=B_2^*u_t(t),\quad t>0.\label{compat}
\end{eqnarray}

If we denote ${\mathcal U}:=(u, u_t, z)^T,$ then
${\mathcal U}^\prime =(u_t, u_{tt}, z_t)^T$ and ${\mathcal U}$ satisfies

\begin{equation}\label{abstract}
\left\{
\begin{array}{l}
{\mathcal U}^\prime ={\mathcal A} {\mathcal U}\\
{\mathcal U}(0)=(u_0, u_1, f^0(-\tau \cdot ))^T,
\end{array}
\right.
\end{equation}
where the operator ${\mathcal A}$ is defined by
\begin{equation}\label{Operator}
{\mathcal A}\left (
\begin{array}{l}
u\\v\\z
\end{array}
\right )
:=\left (
\begin{array}{l}
v\\
-A u-B_1B_1^*v-B_2z(1)\\
-\tau^{-1}z_{\rho}
\end{array}
\right )\,,
\end{equation}
with domain

\begin{equation}\label{dominioOp}
{\mathcal D}({\mathcal A}):=\left\{
\ (u,v,z)^T\in V\times V\times H^1((0,1); U_2)\ :\
A u+B_1B_1^*v\in H \hbox{ and }
 z(0)=B_2^* v
\right\}.
\end{equation}

Denote by ${\mathcal H}$ the Hilbert space
\begin{equation}\label{Hilbert}
{\mathcal H}:=V\times H\times L^2((0,1); U_2),
\end{equation}

equipped with the inner product

\begin{equation}\label{inner}
\left\langle
\left (
\begin{array}{l}
u\\
v\\
z
\end{array}
\right ),\left (
\begin{array}{l}
\tilde u\\
\tilde v\\
\tilde z
\end{array}
\right )
\right\rangle_{\mathcal H}
:=
\langle A^{\frac 1 2} u, A^{\frac 1 2}\tilde u\rangle_H
+\langle v, \tilde v\rangle_H
+\xi\int_0^1\langle z(\rho), \tilde z(\rho)\rangle_{U_2} d\rho ,
\end{equation}
where $\xi$ is any fixed positive number.

The following well--posedness result holds.

\begin{Proposition}\label{WP}
For any initial datum $U_0\in {\mathcal H}$ there exists a unique solution
$U\in C([0,+\infty), {\mathcal H})$  of problem $(\ref{abstract}).$
Moreover, if $U_0\in {\mathcal D}({\mathcal A}),$ then
$$U\in C([0,+\infty), {\mathcal D}({\mathcal A}))\cap C^1([0,+\infty), {\mathcal H}).$$
\end{Proposition}

\noindent {\bf Proof.}
We will show that the operator ${\mathcal A}$ defined by (\ref{Operator}), (\ref{dominioOp}) generates a strongly continuous semigroup   in the Hilbert ${\mathcal H}$ defined in (\ref{Hilbert}), (\ref{inner}).

Denoting by $I$ the identity operator,
we first show that there exists a positive constant $c$ such that ${\mathcal A}-cI$ is dissipative (cfr. \cite{ANP}).
Let $(u,v,z)^T\in {\mathcal D}({\mathcal A}),$  then
$$
\begin{array}{l}
\displaystyle{
\left\langle
{\mathcal A}\left (
\begin{array}{l}
u\\
v\\
z
\end{array}
\right ),\left (
\begin{array}{l}
u\\
v\\
z
\end{array}
\right )
\right\rangle_{\mathcal H}
:=
\left\langle
\left (
\begin{array}{l}
v\\
-A u-B_1B_1^*v-B_2z(1)\\
-\tau^{-1}z_{\rho}
\end{array}
\right ),\left (
\begin{array}{l}
u\\
 v\\
 z
\end{array}
\right )
\right\rangle_{\mathcal H}}
\\
\displaystyle{\hspace{1 cm}
=\langle A^{\frac 1 2}v, A^{\frac 1 2} u\rangle_H
-\langle A u +B_1B_1^*v+B_2 z(1), v\rangle_H
-\xi\int_0^1\langle z_{\rho}(\rho), z(\rho)\rangle_{U_2} d\rho\,.
}\end{array}
$$
Since $A u +B_1B_1^*v+B_2 z(1)\in H\subset V',$ by duality we have
$$
\begin{array}{l}
\displaystyle{
\left\langle
{\mathcal A}\left (
\begin{array}{l}
u\\
v\\
z
\end{array}
\right ),\left (
\begin{array}{l}
u\\
v\\
z
\end{array}
\right )
\right\rangle_{\mathcal H}
= \langle A^{\frac 1 2}v, A^{\frac 1 2} u\rangle_H
-\langle A u, v\rangle_{V^\prime, V} -\langle B_1B_1^*v, v\rangle_{V^\prime, V}
}\\
\hspace{3.5 cm}
\displaystyle{
-\langle B_2 z(1), v\rangle_{V^\prime, V}
-\xi\int_0^1\langle z_{\rho}(\rho), z(\rho)\rangle_{U_2} d\rho
 }\\
\hspace{3 cm} = \displaystyle{
-\Vert B_1^*v\Vert_{U_1}-\langle z(1), B_2^*v\rangle_{U_2}
-\xi\int_0^1\langle z_{\rho}(\rho), z(\rho)\rangle_{U_2} d\rho
\,.
}\end{array}
$$

Integrating by parts and using the relation $z(0)=B_2^* v$, we get
$$\int_0^1\langle z_{\rho}(\rho), z(\rho)\rangle_{U_2} d\rho
=\frac 1 2  (\Vert z(1)\Vert_{U_2}^2-\Vert B_2^*v\Vert_{U_2}^2 ),$$
thus, using also Young's inequality
$$
\left\langle
{\mathcal A}\left (
\begin{array}{l}
u\\
v\\
z
\end{array}
\right ),\left (
\begin{array}{l}
u\\
v\\
z
\end{array}
\right )
\right\rangle_{\mathcal H}
\le \left ( \frac {\xi } 2 +\frac 1 {2\xi } \right )
\Vert B_2^*v\Vert_{U_2}^2\le c\Vert v\Vert_H^2\,,
$$
for a suitable constant $c>0.$
Hence, the operator ${\mathcal A}-cI$ is dissipative.

Now, we show that $\lambda I-{\mathcal A}$ is surjective for some $\lambda >0.$
Given $(f,g,h)^T\in {\mathcal H}$ we seek $(u,v,z)^T\in {\mathcal D}({\mathcal A})$ such that
$$(\lambda I-{\mathcal A})\left (
\begin{array}{l}
u\\
v\\
z
\end{array}
\right )=
\left (
\begin{array}{l}
f\\
g\\
h
\end{array}
\right )\,.$$
This is equivalent to

\begin{eqnarray}
& &\lambda u-v=f\,;\label{S.1}\\
& &\lambda v +A u+B_1B_1^*v +B_2 z(1)=g\,;\label{S.2}\\
& &\lambda z+\tau^{-1}z_{\rho}=h\,.\label{S.3}
\end{eqnarray}

Analogously to \cite{NPSicon06},
suppose that we have found $u$ with the appropriate regularity.
Then, by (\ref{S.1}),
\begin{equation}\label{trovov}
v=\lambda u-f\in V\,.
\end{equation}

Moreover, from (\ref{compat}), (\ref{S.3}) and (\ref{trovov}), $z$ is given by
\begin{equation}\label{trovoz}
z(\rho)=\lambda B_2^*u e^{-\lambda\tau\rho}- B_2^*f e^{-\lambda\tau\rho}+\tau
e^{-\lambda\tau\rho} \int_0^{\rho}  e^{\lambda\sigma\tau} h(\sigma )d\sigma\,,\quad \rho\in (0,1)\,.
\end{equation}
In particular,
\begin{equation}\label{trovoz1}
z(1)=\lambda B_2^*u e^{-\lambda\tau}+z_0\,,
\end{equation}
where
$$z_0=-B_2^*f e^{-\lambda\tau}+\tau e^{-\lambda\tau}\int_0^1  e^{\lambda\sigma\tau}
h(\sigma ) d\sigma\,,
$$
is a fixed element of $U_2$ depending only on $f$ and $h.$

It remains only to determine $u.$
From (\ref{S.1}) and (\ref{S.2}) $u$ satisfies
$$
\lambda^2u+A u+\lambda B_1B_1^* u+B_2z(1)
= g+\lambda f +B_1B_1^*f\,,
$$
and then, by (\ref{trovoz1}),

\begin{equation}\label{trovou}
\lambda^2u+A u+\lambda B_1B_1^* u+\lambda e^{-\lambda\tau }B_2B_2^*u
= g+\lambda f +B_1B_1^*f-B_2 z^0\,.
\end{equation}

We denote the right--hand side of (\ref{trovou}) by $w$, namely
$$w:=g+\lambda f +B_1B_1^*f-B_2 z^0\in H\subset V^\prime\,.$$
Then, from (\ref{trovou}), we have
$$\langle
\lambda^2u+A u+\lambda B_1B_1^* u+\lambda e^{-\lambda\tau }B_2B_2^*u , \varphi \rangle_{V^\prime , V}
=\langle w,\varphi\rangle_{V^\prime , V}.$$

Since $u\in V\subset H,$ we can rewrite
$$
\begin{array}{l}
\langle
\lambda^2u+A u+\lambda B_1B_1^* u+\lambda e^{-\lambda\tau }B_2B_2^*u , \varphi \rangle_{V^\prime , V}\\
\hspace{1 cm} =\lambda^2\langle u, \varphi\rangle_{V^\prime, V}+
\langle A u, \varphi\rangle_{V^\prime, V}+\lambda \langle B_1^*u, B_1^* \varphi\rangle_{V^\prime, V}+\lambda e^{-\lambda\tau }\langle B_2^* u, B_2^* \varphi\rangle_{V^\prime, V}\\
\hspace{1 cm}= \lambda^2\langle u, \varphi\rangle_H +
\langle A^{\frac 1 2} u, A^{\frac 1 2}\varphi\rangle_H+\lambda \langle B_1^*u, B_1^* \varphi\rangle_{U_1}+\lambda e^{-\lambda\tau }\langle B_2^* u, B_2^* \varphi\rangle_{U_2}.
\end{array}
$$
Therefore, we obtain
\begin{equation}\label{var1}
\lambda^2\langle u, \varphi\rangle_H +
\langle A^{\frac 1 2} u, A^{\frac 1 2}\varphi\rangle_H+\lambda \langle B_1^*u, B_1^* \varphi\rangle_{U_1}+\lambda e^{-\lambda\tau }\langle B_2^* u, B_2^* \varphi\rangle_{U_2}=\langle w, \varphi \rangle_{V^\prime , V}.
\end{equation}

The left--hand side of (\ref{var1}) is a continuous and coercive bilinear form on $V.$ Then, Lax--Milgram's lemma implies the existence of a unique solution $u\in V$ of (\ref{var1}) that satisfies
\[
\lambda^2u+A u+\lambda B_1B_1^* u+\lambda e^{-\lambda\tau }B_2B_2^*u 
=w \hbox{ in } V'.
\]
This implies that $Au\in H$ and by defining $v$ by (\ref{trovov}) and $z$ by (\ref{trovoz}), we have found $(u,v,z)^T\in
{\mathcal D}({\mathcal A})$ satisfying (\ref{S.1})--(\ref{S.3}).
This implies that  $\lambda I- {\mathcal A}$ is surjective for all $\lambda >0$ and the same
holds for the operator $\lambda I- ({\mathcal A}-cI).$

Then, the Lumer--Phillips Theorem implies that ${\mathcal A}-cI$ generates a strongly continuous semigroup of contraction in  ${\mathcal H}.$
Hence, the operator  ${\mathcal A}$ generates a  strongly continuous semigroup  in  ${\mathcal H}. \hspace{2 cm}\qed$

\section{Stability result\label{st}}

\hspace{5mm}

\setcounter{equation}{0}


For a fixed constant $\xi$ satisfying
\begin{equation}\label{suxi}
\xi >1,
\end{equation}
we define the energy functional for solutions to problem
(\ref{1.1})--(\ref{1.3}) as
\begin{equation}\label{energy}
E(t):=E(u,t)=\frac 1 2 (\Vert A^{\frac 1 2} u(t)\Vert_H^2+\Vert u_t(t)\Vert_H^2)+
\frac {\xi} 2\int_{t-\tau }^t\Vert B_2^*u_t(s)\Vert_{U_2}^2 ds\,.
\end{equation}
We can obtain a first estimate.

\begin{Proposition}\label{derivE}
For any regular solution of problem $(\ref{1.1})-(\ref{1.3})$
\begin{equation}\label{energy1}
E^\prime (t)\le -\Vert B_1^*u_t(t)\Vert_{U_1}^2+\frac {1+\xi } 2 \Vert B_2^*u_t(t)\Vert_{U_2}^2
-\frac {\xi -1 } 2 \Vert B_2^*u_t(t-\tau)\Vert_{U_2}^2.
\end{equation}
\end{Proposition}

\noindent{\bf Proof:} Differentiating $E(t)$ we get
\begin{eqnarray*}
E'(t)=\langle A^{\frac 1 2} u(t),A^{\frac 1 2} u_t(t)\rangle _H+\langle u_{t}(t), u_{tt}(t)\rangle_H+
\frac{\xi}{2}\Vert B^*_2u_t(t)\Vert_{U_2}^2
-\frac{\xi}{2}\Vert B^*_2u_t(t-\tau)\Vert_{U_2}^2.
\end{eqnarray*}
Hence using the definition of $A$ and (\ref{1.1}) we get successively
\begin{eqnarray*}
E'(t)&=&\langle A u(t), u_{t}(t)\rangle_{V^\prime , V}
-\langle u_t(t), A u(t)\rangle_{V, V^\prime }
-\langle u_t(t), B_1 B_1^*u_{t}(t)\rangle_{V, V^\prime }
\\
& &-\langle u_t(t), B_2B_2^*u_t(t-\tau )\rangle_{V,V^\prime}+
\frac{\xi}{2}\Vert B^*_2u_t(t)\Vert_{U_2}^2
-\frac{\xi}{2}\Vert B^*_2u_t(t-\tau)\Vert_{U_2}^2.
\end{eqnarray*}

Then,
\begin{eqnarray*}
E'(t)=
-\Vert B_1^*u_t(t)\Vert^2_{U_1} -\langle B_2^*u_t (t),B_2^*u_t(t-\tau)\rangle_{U_2}
+\frac{\xi}{2}\Vert B^*_2u_t(t)\Vert_{U_2}^2
-\frac{\xi}{2}\Vert B^*_2u_t(t-\tau)\Vert_{U_2}^2,
\end{eqnarray*}
and (\ref{energy1}) follows from Cauchy-Schwarz's inequality.\hspace{2 cm}\qed

Note that, from (\ref{energy1}), the energy of solutions to problem
(\ref{1.1})--(\ref{1.3}) is not decreasing in general. Indeed the second term
in the right--hand side of (\ref{energy1}), coming from the delay term in (\ref{1.1}), is non negative.
We now consider, as in \cite{pignotti}, the next auxiliary problem which is {\sl close}
to the first one but whose energy is decreasing.

\begin{eqnarray}
& &\varphi_{tt}(t) +A\varphi (t)  + B_1B_1^*\varphi_t(t)+B_2B_2^*\varphi_t(t-\tau )
+\xi B_2B_2^*\varphi_t(t)= 0 ,
\ \quad  t>0,\label{P.1}\\
& &\varphi (0) = \varphi_0,  \quad  \varphi_t(0)=\varphi_1, \quad \label{P.2}\\
& &B_2^*\varphi_t(t)=g^0(t),\quad  t\in (-\tau,0), \label{P.3}
\end{eqnarray}
where $\xi$ is the same constant as in (\ref{energy}).

The well--posedness of system (\ref{P.1})--(\ref{P.3}) can be proved using standard semigroup theory as in Proposition \ref{WP}. Analogously to above we introduce the function
$$\eta (\rho, t)=B_2^*\varphi_t( t-\tau\rho ), \quad \ \rho\in (0,1),\ t>0;$$
and
we rewrite the problem in the abstract form
\begin{equation}\label{abstract2}
\left\{
\begin{array}{l}
\Phi^{\prime}={\mathcal A}^0 \Phi\,,\\
\Phi (0)=(\varphi_0, \varphi_1, g^0(-\tau \cdot ))^T\,,
\end{array}
\right.
\end{equation}
where the operator ${\mathcal A}^0$ is defined by
$$
{\mathcal A}^0\left (
\begin{array}{l}
\varphi\\\psi\\\eta
\end{array}
\right )
:=\left (
\begin{array}{l}
\psi\\
-A \varphi-B_1B_1^*\psi- B_2\eta (1)-\xi B_2B_2^*\psi\\
-\tau^{-1}\eta_{\rho}
\end{array}
\right )\,,$$
with domain ${\mathcal D}({\mathcal A}^0)={\mathcal D}({\mathcal A})$ (see (\ref{dominioOp})) in the Hilbert space ${\mathcal H}$ defined by (\ref{Hilbert}) and (\ref{inner}).

\begin{Proposition}\label{WP2}
For any initial datum $\Phi_0\in {\mathcal H}$ there exists a unique solution
$\Phi\in C([0,+\infty), {\mathcal H})$  of problem $(\ref{abstract2}).$
Moreover, if $\Phi_0\in {\mathcal D}({\mathcal A}^0),$ then
$$\Phi\in C([0,+\infty), {\mathcal D}({\mathcal A}^0))\cap C^1([0,+\infty), {\mathcal H}).$$
\end{Proposition}

For solutions of  problem (\ref{P.1})--(\ref{P.3})  the energy $F(\cdot),$
\begin{equation}\label{energy2}
F(t):=F(\varphi ,t)=
\frac 1 2 (\Vert A^{\frac 1 2} \varphi (t)\Vert_H^2+\Vert \varphi_t(t)\Vert_H^2)+
\frac {\xi} 2\int_{t-\tau }^t\Vert B_2^*\varphi_t(s)\Vert_{U_2}^2 ds\,,
\end{equation}
with $\xi$ satisfying (\ref{suxi}), is decreasing in time.

More precisely, we have the following result.
\begin{Proposition}\label{Problempert}
For any regular solution of problem
$(\ref{P.1})-(\ref{P.3}),$ we have
\begin{equation}\label{stima2}
F^{\prime}(t)\le -\Vert B_1^*\varphi_t(t)\Vert_{U_1}^2
-\frac {\xi -1} 2\Vert B_2^*\varphi_t(t)\Vert_{U_2}^2
-\frac {\xi -1} 2\Vert B_2^*\varphi_t(t-\tau )\Vert_{U_2}^2.
\end{equation}
Then, if $\xi$ satisfies $(\ref{suxi}),$ the energy $F(\cdot)$ is decreasing.
\end{Proposition}

\noindent{\bf Proof.} In order to have (\ref{stima2}) we differentiate (\ref{energy2}).
Hence, using the definition of $A$ and (\ref{P.1}), we obtain
\begin{eqnarray*}
F'(t)&=&\langle A^{\frac 1 2} \varphi ,A^{\frac 1 2} \varphi_t\rangle _H+\langle \varphi_{t}, \varphi_{tt}\rangle_H+
\frac{\xi}{2}\Vert B^*_2\varphi_t(t)\Vert_{U_2}^2
-\frac{\xi}{2}\Vert B^*_2\varphi_t(t-\tau)\Vert_{U_2}^2\\
&=&\langle A \varphi (t), \varphi_{t}(t)\rangle_{V^\prime , V}
-\langle \varphi_t(t), A \varphi (t)\rangle_{V, V^\prime }
-\langle \varphi_t(t), B_1 B_1^*\varphi_{t}(t)\rangle_{V, V^\prime }
- \xi\langle \varphi_t(t), B_2B_2^*\varphi_t(t)\rangle_{V,V^\prime}
\\
& &-\langle \varphi_t(t), B_2B_2^*\varphi_t(t-\tau )\rangle_{V,V^\prime}+
\frac{\xi}{2}\Vert B^*_2\varphi_t(t)\Vert_{U_2}^2
-\frac{\xi}{2}\Vert B^*_2\varphi_t(t-\tau)\Vert_{U_2}^2.
\end{eqnarray*}

Then,
\begin{eqnarray*}
F'(t)&=&
-\Vert B_1^*\varphi_t(t)\Vert^2_{U_1} -\xi \Vert B_2^*\varphi_t(t)\Vert^2_{U_2}
-\langle B_2^*\varphi_t (t),B_2^*\varphi_t(t-\tau)\rangle_{U_2}\\
& &\ +\frac{\xi}{2}\Vert B^*_2\varphi_t(t)\Vert_{U_2}^2
-\frac{\xi}{2}\Vert B^*_2\varphi_t(t-\tau)\Vert_{U_2}^2,
\end{eqnarray*}
and therefore (\ref{energy2}) follows from Cauchy-Schwarz's inequality.\hspace{2 cm}\qed

Consider now the following damped system associated with (\ref{1.1}) and (\ref{1.2}),
 \begin{eqnarray}
& &w_{tt}(t) +A w (t)+B_1 B_1^* w_t=0\quad t>0\label{cons1abstrait}\\
& &w(0)=w_0\quad \mbox{\rm and}\quad w_t(0)=w_1\quad\label{cons2abstrait}
\end{eqnarray}
with $(w_0,w_1)\in V\times H.$
For our stability result we need that this system is exponentially stable or equivalently that the next    observability inequality holds (see Lemma 3.2 of \cite{SergeRendiconti}).
Namely we assume that there exists a time $\overline T>0$ such that for every time $T>\overline T$ there is     a constant $c,$ depending on $T$ but independent of the initial data, such that
\begin{equation}\label{OCabstrait}
E_S(0)\le c\int_0^{T}\Vert B_1^* w_t(t)\Vert_{U_1}^2  dt,
\end{equation}
for every weak solution of problem $(\ref{cons1abstrait}), (\ref{cons2abstrait})$
with initial data $(w_0,w_1)\in V\times H.$

Here $E_S(\cdot)$ denotes the standard energy for wave type equations, that is
$$E_S(t)=E_S(w,t):=\frac 1 2( \Vert A^{\frac 1 2} w(t)\Vert_H^2+\Vert w_t(t)\Vert_H^2).$$

For shortness let us denote by $C_2$ the norm of  $B_2$ 
\begin{equation}\label{C1C2}
\Vert B_2\Vert=\Vert B_2^*\Vert=C_2\,.
\end{equation}

We can prove an exponential stability result for the perturbed problem
(\ref{P.1})--(\ref{P.3}).

\begin{Theorem}\label{stabilitypert}
Assume that $(\ref{suxi})$ holds and that the observability estimate $(\ref{OCabstrait})$  holds for problem $(\ref{cons1abstrait})-(\ref{cons2abstrait})$.
Then, there are two positive constants $K,\tilde \mu$ such that
\begin{equation}\label{exponentialP}
F(t)\le K e^{-\tilde \mu t} F(0),\quad t>0,
\end{equation}
for any solution of problem $(\ref{P.1})-(\ref{P.3}).$
In particular,
\begin{equation}\label{K}
K=\frac{C_0+1}{C_0},
\end{equation}
\begin{equation}\label{mu}
\tilde \mu =\frac{1}{2T}\ln \frac{C_0+1}{C_0},
\end{equation}
with $T$  any fixed time satisfying $T>\max\ \{\overline T,\tau\},$
 $\overline T$ being  an observability time for $(\ref{OCabstrait}),$ and
\begin{equation}\label{Czero}
C_0=\max\  \left\{
2c,\
\frac {32c TC_2^2+\xi } {\xi -1},\
\frac {32c  C_2^2T\xi^2}{\xi -1}
\right\}\,,
\end{equation}
where $C_2$ is as in $(\ref{C1C2})$ and $c:=c(T)$ is the observability constant
in $(\ref{OCabstrait}).$
\end{Theorem}

\noindent {\bf Proof.} Following a classical argument (see \cite{zuazua}) we can decompose the solution $\varphi$ of
 $(\ref{P.1})-(\ref{P.3})$ as
$$\varphi =w+\tilde w$$
where $w$ is the solution of system
$(\ref{cons1abstrait}), (\ref{cons2abstrait})$ with $w_0=\varphi_0,$ $w_1=\varphi_1;$
while $\tilde w$ solves

\begin{eqnarray}
& &\tilde  w_{tt}(t) +A \tilde  w (t)+B_1B_1^* \tilde  w_t(t)=
-\xi B_2B_2^*\varphi_t(t)-B_2B_2^*\varphi_t(t-\tau )
\quad t>0\label{N.1abstrait}\\
& &\tilde  w(0)=0\quad \mbox{\rm and}\quad \tilde   w_t(0)=0\quad\label{N.2abstrait}
\end{eqnarray}

By (\ref{energy2}),

$$
\begin{array}{l}
\displaystyle{
F(0)=E_S(w,0)+\frac{\xi}{2} \int_{-\tau}^0\Vert B_2^*\varphi_t(s)\Vert^2_{U_2}ds}\\\medskip
\hspace{0.9 cm}\displaystyle{ = E_S(w,0)+\frac{\xi}{2} \int_0^{\tau}\Vert B_2^*\varphi_t(t-\tau )\Vert^2_{U_2} dt\,.
}\end{array}
$$
Therefore, from (\ref{OCabstrait}), if $T>\max \{\overline T,\tau\}$ we obtain
\begin{equation}\label{mille}
\begin{array}{l}
\displaystyle{
F(0)\le c\int_0^T\Vert B_1^*w_t(t)\Vert_{U_1}^2 dt +\frac{\xi}{2}
\int_0^T \Vert B_2^*\varphi_t(t-\tau )\Vert^2_{U_2} dt}\\\medskip
\displaystyle{
\hspace{1 cm}\le 2c\int_0^T  (  \Vert B_1^*\varphi_t(t)\Vert_{U_1}^2 +\Vert B_1^*{\tilde w}_t(t)\Vert_{U_1}^2               )dt +\frac{\xi}{2}
\int_0^T\Vert B_2^*\varphi_t(t-\tau )\Vert^2_{U_2}dt,}
\end{array}
\end{equation}
where $c$ is the observability constant for the damped       system (\ref{cons1abstrait}), (\ref{cons2abstrait}).

Now, observe that from (\ref{N.1abstrait}),
$$
\begin{array}{l}
\displaystyle{
\frac{d}{dt} \frac 1 2 (\Vert \tilde w_t(t)\Vert_H^2+\Vert A^{\frac 1 2}\tilde w(t)\Vert^2_H)+\|B_1^* \tilde w_t\|_{U_1}^2 =
\langle \tilde w_t, \tilde w_{tt}+A\tilde w+B_1B_1^* \tilde  w_t\rangle}\\
\hspace{1 cm}
= \langle \tilde w_t, -
\xi B_2B_2^*\varphi_t(t)-B_2B_2^*\varphi_t(t-\tau )
\rangle_H\,.
\end{array}
$$

Integrating in time from 0 to  $t,$ for $t\in (0,2T],$
and using (\ref{N.2abstrait}) we have
$$
\begin{array}{l}
\displaystyle{
\frac 1 2 (\Vert \tilde w_t(t)\Vert_H^2+\Vert A^{\frac 1 2}\tilde w(t)\Vert^2_H) +\int_0^t\|B_1^* \tilde w_t(s)\|_{U_1}^2\,ds }
\\
\hspace{1 cm}
= \displaystyle{\int_0^t \langle \tilde w_t, -
\xi B_2B_2^*\varphi_t(s)-B_2B_2^*\varphi_t(s-\tau )
\rangle_H ds}\,,
\end{array}
$$
and then
\begin{equation}\label{mille1}
\begin{array}{l}
\displaystyle{
\Vert \tilde w_t(t)\Vert_H^2+2\int_0^t\|B_1^* \tilde w_t(s)\|_{U_1}^2\,ds \le
 \frac 1 {8TC_2^2}\int_0^t \Vert B_2^*\tilde w_t(s)\Vert_{U_2}^2 ds
+8 TC_2^2\xi^2\int_0^t \Vert B_2^* \varphi_t(s)\Vert_{U_2}^2 ds
}\\
\hspace{1 cm}\displaystyle{
+ \frac 1 {8TC_2^2}\int_0^t \Vert B_2^*\tilde w_t(s)\Vert_{U_2}^2 ds
+8 TC_2^2\int_0^t \Vert B_2^* \varphi_t(s-\tau )\Vert_{U_2}^2 ds
}\,,
\end{array}
\end{equation}
where $C_2$ was defined in (\ref{C1C2}).

This estimate directly implies that  for all $t\in [0,2T],$ one has

\begin{equation}\label{mille2}
\begin{array}{l}
\displaystyle{
\Vert \tilde w_t(t)\Vert_H^2+2\int_0^t\|B_1^* \tilde w_t(s)\|_{U_1}^2\,ds \le
 \frac 1 {4TC_2^2}\int_0^{2T} \Vert B_2^*\tilde w_t(s)\Vert_{U_2}^2 ds
+8 TC_2^2\xi^2\int_0^{2T} \Vert B_2^*\varphi_t(s)\Vert_{U_2}^2 ds
}\\
\hspace{1 cm}\displaystyle{
+8 TC_2^2\int_0^{2T} \Vert B_2^* \varphi_t(s-\tau )\Vert_{U_2}^2 ds
}\,,
\end{array}
\end{equation}
and so, integrating in $[0,2T],$
\begin{equation}\label{mille3}
\begin{array}{l}
\displaystyle{\int_0^{2T}
\Vert \tilde w_t(t)\Vert_H^2 dt +2\int_0^{2T}\int_0^t\|B_1^* \tilde w_t(s)\|_{U_1}^2\,ds\,dt\le
 \frac 1 {2C_2^2}\int_0^{2T} \Vert B_2^*\tilde w_t(s)\Vert_{U_2}^2 ds}\\
\hspace{1 cm}\displaystyle{+16 T^2C_2^2\xi^2\int_0^{2T} \Vert B_2^* \varphi_t(s)\Vert_{U_2}^2 ds
+16 T^2C_2^2\int_0^{2T} \Vert B_2^* \varphi_t(s-\tau )\Vert_{U_2}^2 ds
}\,.
\end{array}
\end{equation}

Therefore,
$$
\begin{array}{l}
\displaystyle{\frac12\int_0^{2T}
\Vert \tilde w_t(t)\Vert_H^2 dt +2\int_0^{2T}\int_0^t\|B_1^* \tilde w_t(s)\|_{U_1}^2\,ds\,dt\le
 16 T^2C_2^2\xi^2\int_0^{2T} \Vert B_2^* \varphi_t(s)\Vert_{U_2}^2 ds
}\\
\hspace{1 cm}\displaystyle{
+16 T^2C_2^2\int_0^{2T} \Vert B_2^* \varphi_t(s-\tau )\Vert_{U_2}^2 ds
}\,,
\end{array}
$$
from which follows
\[
\begin{array}{l}
\displaystyle{\int_0^{2T}\int_0^t\|B_1^* \tilde w_t(s)\|_{U_1}^2\,ds\,dt\le
 8 T^2C_2^2\xi^2\int_0^{2T} \Vert B_2^* \varphi_t(s)\Vert_{U_2}^2 ds
}\\
\hspace{1 cm}\displaystyle{
+8 T^2C_2^2\int_0^{2T} \Vert B_2^* \varphi_t(s-\tau )\Vert_{U_2}^2 ds
}\,.
\end{array}
\]
Using the fact that
\begin{eqnarray*}
\int_0^{2T}\int_0^t\|B_1^* \tilde w_t(s)\|_{U_1}^2\,ds\,dt&=&
\int_0^{2T} \|B_1^* \tilde w_t(s)\|_{U_1}^2(2T-s)\,ds\\
&\geq &\int_0^{T} \|B_1^* \tilde w_t(s)\|_{U_1}^2(2T-s)\,ds
\\
&\geq &T \int_0^{T} \|B_1^* \tilde w_t(s)\|_{U_1}^2\,ds,
\end{eqnarray*}
we deduce that
\begin{equation}\label{fe1}
\begin{array}{l}
\displaystyle{\int_0^{T}\|B_1^* \tilde w_t(s)\|_{U_1}^2\,ds\le
 8 TC_2^2\xi^2\int_0^{2T} \Vert B_2^* \varphi_t(s)\Vert_{U_2}^2 ds
}\\
\hspace{1 cm}\displaystyle{
+8 TC_2^2\int_0^{2T} \Vert B_2^* \varphi_t(s-\tau )\Vert_{U_2}^2 ds
}\,.
\end{array}
\end{equation}

Using (\ref{fe1}) in (\ref{mille}) we obtain
\begin{equation}\label{fe2}
\begin{array}{l}
\displaystyle{
F(0)\le 
2c\int_0^T\Vert B_1^*\varphi_t(t)\Vert_{U_1}^2 dt+(16c TC_2^2+\frac {\xi } 2)\int_0^{2T}\Vert B_2^*\varphi (t-\tau )\Vert_{U_2}^2 dt
}\\
\hspace{1.3 cm}\displaystyle{
+16c  C_2^2T\xi^2 \int_0^{2T}\Vert B_2^*\varphi_t(t)\Vert_{U_2}^2 dt
},
\end{array}
\end{equation}
that we rewrite as
\begin{equation}\label{fe3}
\begin{array}{l}
\displaystyle{
F(0)\le 2c\int_0^{2T}\Vert B_1^*\varphi_t(t)\Vert_{U_1}^2 dt
+\frac {32c TC_2^2+\xi } {\xi -1}
\left (
\frac {\xi -1} 2\int_0^{2T}\Vert B_2^*\varphi (t-\tau )\Vert_{U_2}^2 dt\right )
}\\
\hspace{2 cm}\displaystyle{+
\frac {32c  C_2^2T\xi^2}{\xi -1}
\left (
\frac {\xi -1} 2\int_0^{2T}\Vert B_2^*\varphi (t)\Vert_{U_2}^2 dt\right )
\le -C_0\int_0^{2T} F^\prime (t) dt},
\end{array}
\end{equation}
with $C_0$ as in (\ref{Czero}).

Therefore, from (\ref{fe3}), using also  that $F(\cdot)$ is decreasing we obtain
$$F(2T)\le F(0)\le C_0 (F(0)-F(2T)).$$

Then,
$$F(2T)\le \frac {C_0}{C_0+1} F(0),$$
and this implies the exponential estimate (\ref{exponentialP}) with $K, \tilde\mu$ as in (\ref{K}), (\ref{mu}),
due to the semigroup property together with the fact that $F$ is non increasing.\hspace{2 cm}\qed

Now, let us recall the following classical result of Pazy (Theorem 1.1 in Ch. 3 of \cite{pazy}).
\begin{Theorem}\label{Pazy}
Let $X$ be a Banach space and let $A$ be the infinitesimal generator of a $C_0$ semigroup $T(t)$ on $X,$ satisfying $\Vert T(t)\Vert \le Me^{\omega t}.$
If $B$ is a bounded linear operator on $X$ then $A+B$ is the infinitesimal generator of a $C_0$ semigroup $S(t)$ on $X,$ satisfying
$\Vert S(t)\Vert\le M e^{(\omega+M\Vert B\Vert )t}\,.$
\end{Theorem}

Using this perturbation result we are ready to give the asymptotic stability result for problem
(\ref{1.1})--(\ref{1.3}) when the norm of the delay feedback $B_2^*$ is sufficiently small.

\begin{Theorem}\label{exponential}
Assume that the observability estimate $(\ref{OCabstrait})$  holds for problem $(\ref{cons1abstrait})-(\ref{cons2abstrait})$.
For all $\xi>1$ in the definition $(\ref{energy})$, there is $\beta >0$ depending on $\bar T$, $\tau$, $\xi$ and on the operator $B_1,$  such that if the delay feedback satisfies
$\Vert B_2^*\Vert <\beta $, then 
there exist positive constants $K, \mu$ for which  we have
\begin{equation}\label{stimaOr}
E(t)\le K e^{- \mu t} E(0),\quad t>0,
\end{equation}
for any solution of $(\ref{1.1})-(\ref{1.3}).$
\end{Theorem}

\noindent{\bf Proof.}
We can see problem (\ref{1.1})--(\ref{1.3}) as a perturbation of the auxiliary one.
Therefore,
$${\mathcal A}\left (
\begin{array}{l}
u\\v\\z
\end{array}
\right )
= ({\mathcal A}^0 +{\mathcal B})\left (
\begin{array}{l}
u\\v\\z
\end{array}
\right )$$
with

$${\mathcal B}\left (
\begin{array}{l}
u\\v\\z
\end{array}
\right )
=
\left (
\begin{array}{l}
0\\-\xi B_2B_2^*v\\0
\end{array}
\right )\,.$$

From Theorem \ref{stabilitypert} and Theorem \ref{Pazy}, we know that
if
\begin{equation}\label{dieci}
-\tilde \mu +K\Vert {\mathcal B}\Vert<0,
\end{equation}
where $\tilde\mu$ and $K$ are defined by  (\ref{K}),
(\ref{mu}) and (\ref{Czero}),
then
$$E(t)\le K e^{-\mu t} E(0),$$
with $\mu =\tilde\mu-K\Vert {\mathcal B}\Vert,$ for any solution of problem (\ref{1.1})--(\ref{1.3}).

It remains to prove that (\ref{dieci}) is satisfied for $\Vert B_2^*\Vert $ sufficiently small.
We can rewrite (\ref{dieci}) as
$$\xi \Vert B_2\Vert^2 <\frac {\tilde\mu }{K},$$
that is
\begin{equation}\label{oggi1}
\xi C_2^2 <\frac{1}{2 T} \frac {C_0}{C_0+1}\ln \frac {C_0+1}{C_0}\,.
\end{equation}

The difficulty is that the constant $C_0$ (defined by (\ref{Czero})) appearing in the right--hand side of this estimate depends on $\xi$ and $C_2$ as well. 
So let us consider the continuous  function
$h:(0,+\infty)\rightarrow (0,+\infty),$
$$ h(s):= \frac{s}{s+1} \ln \frac{s+1}{s}\,.$$
Then, $h$ tends to zero for $s\rightarrow 0^+$ and for $s\rightarrow +\infty .$
Moreover, $h$ assumes the maximal value $1/e$ at $\frac 1 {e-1}$, is increasing before $\frac 1 {e-1}$ and decreasing after.

Now 
it is easy to check that $\frac \xi {\xi -1}>\frac 1 {e-1}.$
Considering that $\xi$ is fixed $>1$ and $T$ is fixed as well, we consider 
$C_0$ as a function of $C_2\geq 0$, that we write $C_0(C_2)$.
But we remark that from its definition, $C_0$ is non decreasing (in $C_2$) and 
\[
C_0(0)= \max\  \left\{
2c,\
\frac{\xi }{\xi-1}
\right\}>\frac 1 {e-1}.
\]
Hence 
$h(C_0(C_2))$ is non increasing as a function of $C_2$ with $h(C_0(0))>0$ and
since the left--hand side of (\ref{oggi1}) is increasing in $C_2$ and is zero at $C_2=0$, there exists a point $\beta>0$ such that
\[
2\xi \beta^2 T=h(C_0(\beta)),
\]
and 
for which (\ref{oggi1}) holds for all $C_2\in [0,\beta)$.

Obviously $\beta$ depends on $T$ (and then on $\bar T$ and $\tau$), on $\xi$
and, through the observability constant $c$ and the time $T,$
on the feedback operator $B_1.\hspace{2 cm}\qed$
 
\begin{Remark}\label{rkB1borne}
{\rm
If $B_1$ is bounded, namely if $B_1\in \mathcal{L}(U_1,H)$, then by Proposition 1 of \cite{Haraux:89}, the system $(\ref{cons1abstrait})-(\ref{cons2abstrait})$ is exponentially stable (or equivalently the observability estimate $(\ref{OCabstrait})$  holds for problem $(\ref{cons1abstrait})-(\ref{cons2abstrait})$) if and only if the observability estimate
\begin{equation}\label{OCabstraitconservative}
E_S(0)=\frac 1 2( \Vert A^{\frac 1 2} w_0\Vert_H^2+\Vert w_1\Vert_H^2)\le c\int_0^{T}\Vert B_1^* \varphi_t(t)\Vert_{U_1}^2  dt,
\end{equation}
holds for some $T>0$ and $c>0$,
for every weak solution $\varphi$ of the conservative system
\begin{eqnarray}
& &\varphi_{tt}(t) +A \varphi (t)=0\quad t>0\label{conservative1abstrait}\\
& &\varphi(0)=w_0\quad \mbox{\rm and}\quad \varphi_t(0)=w_1\quad\label{conservative2abstrait}
\end{eqnarray}
with initial data $(w_0,w_1)\in V\times H.$
\hspace{2 cm}$\qed$
}
\end{Remark}

\section{The wave equation\label{WW}}

\hspace{5mm}

\setcounter{equation}{0}

\subsection{Internal dampings\label{WWint}}

Our first application concerns the wave equation with locally distributed internal dampings.
More precisely, let $\Omega\subset\RR^n$ be an open bounded domain   with  a Lipschitz boundary
$\partial\Omega$. We suppose given $b_1,b_2$ in $L^\infty (\Omega)$ such that
$$b_1(x), b_2(x)\ge 0 \quad {\mbox {\rm a.e.}}\ x\in \Omega.$$

Let us consider the initial boundary value  problem

 \begin{eqnarray}
& &u_{tt}(x,t) -\Delta u (x,t)+b_1(x) u_t(x,t)+b_2(x) u_t(x,t-\tau)=0\ \ \mbox{\rm in}\quad\Omega\times
(0,+\infty),\label{W.1}\\
& &u (x,t) =0\quad \mbox{\rm on}\quad\partial\Omega\times
(0,+\infty),\label{W.2}\\
& &u(x,0)=u_0(x)\quad \mbox{\rm and}\quad u_t(x,0)=u_1(x)\quad \hbox{\rm
in}\quad\Omega ,\label{W.3}\\
& &\sqrt{b_2} u_t(x,t)=f^0(x,t) \quad \mbox{\rm in}\quad\omega_2\times
(-\tau, 0),\label{Iw}\\
\end{eqnarray}
with
initial
data $(u_0, u_1, f^0)\in H^1_0(\Omega)\times L^2(\Omega)\times L^2((-\tau,0);L^2(\omega_2))$, where $\omega_i=\{x\in \Omega: b_i(x)>0\}$ is the support of $b_i$, $i=1$ or 2.

This problem enters into our previous framework, if we take
$H=L^2(\Omega)$ and the operator  $A$ defined by
$$A:{\mathcal D}(A)\rightarrow H\,:\,  u\rightarrow -\Delta u,$$
where ${\mathcal D}(A)=\{u\in H^1_0(\Omega): \Delta u\in L^2(\Omega)\}.$
This operator $A$ is a self--adjoint and positive operator with a compact inverse in $H$
and is such that $V={\mathcal D}(A^{1/2})=H^1_0(\Omega).$
We then define $U_1=L^2(\omega_1)$, $U_2=L^2(\omega_2)$ and the operators $B_i, i=1,2,$ as
\be\label{defBi}
B_i:U_i\rightarrow H: \quad v\rightarrow \sqrt{b_i(x)}\tilde v ,
\ee
where $\tilde v\in L^2(\Omega)$ is the extension of $v$ by zero outside $\omega_i.$
It is easy to verify that
$$B_i^* \varphi =\sqrt {b_i} \varphi_{\vert_{\omega_i}}\quad\mbox{\rm for}\ \varphi\in H.$$
As
$B_iB_i^* \varphi=b_i\varphi ,$ for any $\varphi\in H$ and $i=1,2,$
we deduce  that problem (\ref{W.1})--(\ref{Iw}) enters in the abstract framework
(\ref{1.1})--(\ref{1.3}).

In this setting, the energy functional is

\begin{equation}\label{ener}
 \displaystyle{
E(t)= \frac{1}{2}\int_\Omega \{ u_t^2(x,t) +\vert \nabla u(x,t)\vert^2\} dx +
 \frac \xi 2\int_{t-\tau}^{t} \int_{\Omega}b_2(x) u_t^2(x,s) dx ds,
}
\end{equation}
which is the standard energy for wave equation
$$E_S(t)=E_S(w,t):=\frac 1 2 \int_{\Omega}( w_t^2 +\vert\nabla w\vert^2) dx,$$
plus an integral term
due to the presence of a time delay.

Since $B_1$ is bounded, according to Remark \ref{rkB1borne}, our main assumption concerns the existence of an observability estimate for the standard wave equation:

 \begin{eqnarray}
& &\varphi_{tt}(x,t) -\Delta \varphi (x,t)=0\quad \mbox{\rm in}\quad\Omega\times
(0,+\infty)\label{C.1}\\
& &\varphi (x,t) =0\quad \mbox{\rm on}\quad\partial\Omega\times
(0,+\infty)\label{C.2}\\
& &\varphi(x,0)=w_0(x)\quad \mbox{\rm and}\quad \varphi_t(x,0)=w_1(x)\quad \hbox{\rm
in}\quad\Omega\label{C.3}
\end{eqnarray}
with $(w_0,w_1)\in H_0^1(\Omega)\times L^2(\Omega).$

We then assume that
there exists a time $\overline T>0$ such that for every time $T>\overline T$
there is a
constant $c,$ depending on $T$ but independent of the initial data, such that
\begin{equation}\label{OC}
E_S(0)\le c\int_0^{T}\int_{\Omega}b_1(x)\varphi_t^2(x,s) dx ds,
\end{equation}
for every weak solution of problem $(\ref{C.1})-(\ref{C.3}).$

According to (\ref{C1C2}) we have
\begin{equation}\label{C1C2wave}
\|B_2\|=\Vert b_2\Vert_{\infty}^{1/2},
\end{equation}
where, for $v\in L^\infty (\Omega),$ we denote 
 $\Vert v\Vert_{\infty}=\sup_{x\in \Omega} |v(x)|$,
the $L^\infty$ norm of  $v.$

Therefore, according to Theorem \ref{exponential}, we have the next result:

\begin{Theorem}\label{exponentialw}
Assume that the observability estimate $(\ref{OC})$  holds for the wave equation $(\ref{C.1})-(\ref{C.3}).$ For all $\xi>1$ in the definition $(\ref{ener})$, there is $\beta >0$ depending on $\bar T$, $\tau$, $\xi$ and $b_1$  such that if 
$\Vert b_2\Vert_{\infty} <\beta $, then 
there exist positive constants $K, \mu$ for which  we have
$$
E(t)\le K e^{- \mu t} E(0),\quad t>0,
$$
for any solution of $(\ref{W.1})-(\ref{Iw}).$
\end{Theorem}

\begin{Remark}\label{rkwave1}
{\rm
1. From Lemma VII.2.4 of  \cite{Lions} (see also \cite{Komornikbook, Lag83, LT, liu, zuazua}), the 
observability estimate (\ref{OC})  holds for the wave equation $(\ref{C.1})-(\ref{C.3})$
if    the boundary of $\Omega$ is
of class $C^2$, if $T$ is bigger than the diameter of $\Omega$
and if 
\begin{equation}\label{perobs}
b_1(x)\ge b^0>0,\quad {\mbox {\rm a.e.}}\  x\in \omega,
\end{equation}
when the open subset
$\omega$ of $\Omega$ is   a   neighborhood of $\bar \Gamma_0$, where
\be\label{defgamma0}\Gamma_0=\{\, x\in\partial\Omega\, :\ (x-x_0)\cdot \nu (x)>0\,\},
\ee
for some $x_0\in {\mathbb R}^n$ and $\nu(x)$ is the outer unit normal vector at $x\in\partial\Omega.$
\\
2. From \cite{BLR}, the 
observability estimate (\ref{OC})  also holds for the wave equation $(\ref{C.1})-(\ref{C.3})$
if    the boundary of $\Omega$ is
of class $C^\infty$ and if (\ref{perobs}) holds when the open subset
$\omega$ of $\Omega$ satisfies the geometric control property.
\hspace{2 cm}$\qed$
}
\end{Remark}

\begin{Remark}
{\rm
According to point 1 of the previous remark, Theorem \ref{exponentialw} allows to recover the results from Theorem 1.2 of \cite{pignotti} in a larger setting.
\hspace{2 cm}$\qed$
}
\end{Remark}

\subsection{Internal and boundary dampings \label{Wavebound}}

We assume here that the boundary 
$\partial\Omega$ of $\Omega$ is splitted up as 
$\partial \Omega =  \Gamma_0 \cup {\Gamma}_1,$
where $\Gamma_0,\ \Gamma_1$ are closed subsets of $\partial \Omega$
with   $\Gamma_0\cap
\Gamma_1=\emptyset$. Moreover we assume that   $\Gamma_0$ and  $\Gamma_1$ have an non empty interior (on $\partial \Omega$). We suppose given 
$k\in L^\infty (\Gamma_0)$ and   $b \in L^\infty (\Omega)$ such that
 $b(x)\ge 0$ a.e. $x\in\Omega$ and 
$$
k(x)\ge k_0> 0 \hbox{ a.e. } x\in\Gamma_0.$$

We here consider 
the problem
\begin{eqnarray}
& &u_{tt}(x,t)-\Delta u(x,t)  + b(x) u_t(x,t-\tau)= 0, \quad x\in\Omega,\ t>0,\label{a1}\\
& &u(x,t)=0, \quad x\in  \Gamma_1,\ t> 0\label{a2}\\
& &\frac{\partial u}{\partial \nu}(x,t)=- k(x) u_t(x,t), \quad x\in  \Gamma_0,\ t> 0\label{a3}\\
& &u(x,0) = u_0(x),  \  u_t(x,0)=u_1(x),\quad  x\in
\Omega,\label{a4}\\
& &\sqrt{b} u_t(x,t)=f^0(x,t),\quad   x\in\Omega, \ t\in (-\tau,0), \label{a5}
\end{eqnarray}
with initial data 
in a suitable space.

This problem enters into our previous framework, if we take
$H=L^2(\Omega)$ and the operator  $A$ defined by
$$A:{\mathcal D}(A)\rightarrow H\,:\,  u\rightarrow -\Delta u,$$
where 
$${\mathcal D}(A):=\{\, u\in H^1_{\Gamma_1}(\Omega) 
\,:\, \Delta u\in L^2(\Omega)\hbox{ and } \frac {\partial u}{\partial\nu}=0\ \mbox{\rm on}\ \Gamma_0\,\}
,$$
with
$$H_{\Gamma_1}^1:=\{\, u\in H^1(\Omega)\,:\, u=0\ \mbox{\rm on}\ \Gamma_1\,\}.$$

We then define $U_1:=L^2(\Gamma_0),$ $U_2:=L^2(\omega_2)$ ($\omega_2$ being the support of $b$)
and the operators $B_1, B_2$ as
$$B_2 \in {\cal L}(U_2; H), \, B_2u  = \sqrt{b(x)} \, \tilde u,
\  \forall \, u \in L^2(\omega_2),$$
and
$$B_1 \in {\cal L}(U_1;V^\prime    ), \, B_1u = \sqrt{k} \, A_{-1} Nu, \, \forall \, u \in L^2(\Gamma_0), \,
B_1^*w =\sqrt{k}  w_{|\Gamma_0}, \, \forall \, w \in 
V:=
{\cal D}(A^{1/2}),$$          where                                                                      $A_{-1}$ is the extension of $A$ to $H$,
namely for all $h\in H$ and $\varphi \in {\cal D}(A)$, $A_{-1}h$ is
the unique element in $({\mathcal D}(A))^\prime$ 
                                (the duality is
in the sense of $H$),
      such that (see for instance
\cite{TucsnakWeiss})
\[
\langle A_{-1}h; \varphi\rangle_{({\mathcal D}(A))^\prime,  {\mathcal D}(A)}                  =\int_\Omega h
A\varphi\, dx.\]

Here and below $N \in {\cal L}(L^2(\Gamma_0);L^2(\Omega))$ is defined as follows:  for all
$v \in L^2(\Gamma_0), \, Nv$ is the unique solution (transposition
solution) of
$$\Delta Nv =0, \,
Nv_{|\Gamma_1} = 0, \,
\frac{\partial Nv}{\partial \nu}_{|\Gamma_0} = v.$$

With these definitions, we can show that problem (\ref{a1})--(\ref{a5}) enters in the abstract framework
(\ref{1.1})--(\ref{1.3}).

Now, the energy functional is

\begin{equation}\label{enerWB}
 \displaystyle{
E(t)= \frac{1}{2}\int_\Omega \{ u_t^2(x,t) +\vert \nabla u(x,t)\vert^2\} dx +
 \frac \xi 2\int_{t-\tau}^{t} \int_{\Omega}b(x) u_t^2(x,s) dx ds.
}
\end{equation}

As $B_1$ is not bounded, we need to consider the non delayed  system

 \begin{eqnarray}
& &w_{tt}(x,t) -\Delta w (x,t)=0\quad \mbox{\rm in}\quad\Omega\times
(0,+\infty)\label{C.1B}\\
& &w (x,t) =0\quad \mbox{\rm on}\quad\Gamma_1\times
(0,+\infty)\label{C.2B}\\
& &\frac{\partial w}{\partial \nu}(x,t)=- k(x) w_t(x,t), \quad x\in  \Gamma_0,\ t> 0\label{a3B}\\
& &w(x,0)=w_0(x)\quad \mbox{\rm and}\quad w_t(x,0)=w_1(x)\quad \hbox{\rm
in}\quad\Omega\label{C.3B}
\end{eqnarray}
with $(w_0,w_1)\in H_{\Gamma_1}^1(\Omega)\times L^2(\Omega).$

Hence our main  assumption will be:
There exists a time $\overline T>0$ such that for every time $T>\overline T$
there is a
constant $c,$ depending on $T$ but independent of the initial data, such that
\begin{equation}\label{OCB}
E_S(0)\le c\int_0^{T}\int_{\Gamma_0}k(x) w_t^2(x,s) dx ds,
\end{equation}
for every weak solution of problem $(\ref{C.1B})-(\ref{C.3B}).$

Then, our previous results apply also to this model and we can restate 
Theorem \ref{exponential}.
\begin{Theorem}\label{exponentialwB}
Assume that the observability estimate $(\ref{OCB})$ holds
for every weak solution of problem $(\ref{C.1B})-(\ref{C.3B}).$
For all $\xi>1$ in the definition $(\ref{enerWB})$, there is $\beta >0$ depending on $\bar T$, $\tau$, $\xi$ and $k,$  such that if 
$\Vert b(x)\Vert_{\infty} <\beta $, then 
there exist positive constants $K, \mu$ for which  we have
$$
E(t)\le K e^{- \mu t} E(0),\quad t>0,
$$
for any solution of $(\ref{a1})-(\ref{a5}).$
\end{Theorem}

\begin{Remark}\label{rkwavebdy}
{\rm
1. From Theorem 1 and Remark 1 of  \cite{komornikzuazua} (see also \cite{Chen79,Chen81,Chen81b,LaJDE83}), the 
observability estimate (\ref{OCB})  holds for the damped wave equation $(\ref{C.1B})-(\ref{C.3B})$
if    the boundary of $\Omega$ is
of class $C^2$, if $T$ is large enough
and if $\Gamma_0$ is given by (\ref{defgamma0})
for some $x_0\in {\mathbb R}^n$.
\\
2. From \cite{BLR}, the 
observability estimate (\ref{OCB})  also holds for the damped wave equation $(\ref{C.1B})-(\ref{C.3B})$
if the boundary of $\Omega$ is
of class $C^\infty$ and if the part $\Gamma_0$ satisfies  the geometric control property.
\\
3. If we suppress the assumption $\Gamma_0\cap
\Gamma_1=\emptyset$, then Theorem 1  of  \cite{komornikzuazua}  shows that
the 
observability estimate (\ref{OCB})  holds for the damped wave equation $(\ref{C.1B})-(\ref{C.3B})$
under the same assumptions than in point 1 but with the choice
$k(x)=(x-x_0)\cdot \nu(x)$ and if $n\leq 3$ (see also Proposition 6.4 of \cite{grisvard:89b} in dimension 2). For this example,  $k$ is no more   uniformly positive  on $\Gamma_0$, nevertheless it enters into our abstract framework.
}\hspace{2 cm}\qed
\end{Remark}

\begin{Remark}\label{analogo}
{\rm
This result, namely exponential decay of the energy for solutions
   to problem (\ref{a1})--(\ref{a5})
for ``small'' internal delay feedback,
      has been first proved in \cite{ANP},   for $b$ and $k$ constant and $\Gamma_0$  given by (\ref{defgamma0}), by 
constructing
a suitable Lyapunov functional and using the multiplier method.
We  give here a simpler proof by using a more general method, allowing to weaken the assumptions on $b$, $k$ and $\Gamma_0$.
\hspace{2 cm}$\qed$
}
\end{Remark}
\section{The elasticity system\label{Elastic}}

\hspace{5mm}

\setcounter{equation}{0}

\subsection{Internal dampings \label{Elasticint}}

Here we consider the following elastodynamic system

 \begin{eqnarray}
\nonumber
& &u_{tt}(x,t) -\mu \Delta  u (x,t)-(\lambda+\mu) \nabla \hbox{ div }u
\\
&&\hspace{3cm}+b_1(x)  u_t(x,t)+b_2(x) 
u_t(x,t-\tau)=0\quad \mbox{\rm in}\quad\Omega\times
(0,+\infty),\label{P.1ela}\\
& &u (x,t) =0\quad \mbox{\rm on}\quad\partial\Omega\times
(0,+\infty),\label{P.2ela}\\
& &u(x,0)=u_0(x)\quad \mbox{\rm and}\quad u_t(x,0)=u_1(x)\quad \hbox{\rm
in}\quad\Omega ,\label{P.3ela}\\
& &\sqrt{b_2} u_t(x,t)=f^0(t)\quad \mbox{\rm in} \ \omega_2\times (-\tau , 0), \label{Ie}
\end{eqnarray}
with
initial
data $(u_0, u_1, f^0)\in H^1_0(\Omega)^n\times L^2(\Omega)^n\times L^2((-\tau ,0);L^2(\omega_2)^n)$
and
$b_1,b_2$ satisfying
the
  same assumptions as in subsection \ref{WWint}. Note that in this case the state variable
$u$ is vector-valued and $\lambda,\mu$ are the Lam\'e coefficients that are positive real numbers.

As before this problem enters into our abstract setting, once we take
$
H=L^2(\Omega)^n$,
and
   $A$ defined by
$$A:{\mathcal D}(A)\rightarrow H\,:\,  u\rightarrow -\mu \Delta  u (x,t)-(\lambda+\mu) \nabla \hbox{ div }u,$$
where ${\mathcal D}(A)=\{u\in H^1_0(\Omega)^n: \mu \Delta  u +(\lambda+\mu) \nabla \hbox{ div }u
\in L^2(\Omega)^n\}.$

The operator $A$ is a self--adjoint and positive operator with a compact inverse in $H$
and is such that $V={\mathcal D}(A^{1/2})=H^1_0(\Omega)^n$
equipped with the inner product
$$
(u,v)_V=\int_\Omega \Big(\mu\sum_{i,j=1}^n \partial_i u_j \partial_i v_j+(\lambda+\mu)  \hbox{ div }u\hbox{ div }v\Big)\,dx, \ \ \forall u,v\in H^1_0(\Omega)^n.
$$

We then define $U_i=L^2(\omega_i)^n$ and the operators $B_i, i=1,2,$ as
$$B_i:U_i\rightarrow H : \quad v\rightarrow \sqrt{b_i}\tilde v ,$$
where $\tilde v$ is the extension of $v$ by zero outside $\omega_i.$
As before
$$B_i^*(\varphi )=\sqrt {b_i} \varphi_{\vert_{\omega_i}}\quad\mbox{\rm for}\ \varphi\in H,$$
and thus $B_iB_i^*(\varphi)=b_i\varphi,$ for any $\varphi\in H$ and $i=1,2.$
So, problem (\ref{P.1ela})--(\ref{Ie}) enters in the abstract framework
(\ref{1.1})--(\ref{1.3}).

Therefore in order to apply the abstract results of section \ref{st}, we only need to check the
observability estimate for the associated conservative   system: There exists a time $T>0$ and a constant $c>0$ such that
\begin{equation}\label{OEelasticityint}
\frac12((w_0,w_0)_V+\int_\Omega |w_1|^2\,dx)
\leq c \int_0^{T}\int_{\Omega}b_1(x)|\varphi_t|^2(x,s) dx ds,
\end{equation}
for every weak solution $\varphi$ of
\begin{eqnarray*}
& &\varphi_{tt}(x,t) -\mu \Delta  \varphi (x,t)-(\lambda+\mu) \nabla \hbox{ div }\varphi(x,t)=0\quad \mbox{\rm in}\quad\Omega\times
(0,+\infty), \\
& &\varphi (x,t) =0\quad \mbox{\rm on}\quad\partial\Omega\times
(0,+\infty), \\
& &\varphi(x,0)=w_0(x)\quad \mbox{\rm and}\quad \varphi_t(x,0)=w_1(x)\quad \hbox{\rm
in}\quad\Omega,
\end{eqnarray*}
with
initial
data $(w_0, w_1)\in H^1_0(\Omega)^n\times L^2(\Omega)^n$.

If such an estimate holds, the stability result from section \ref{st} can be applied to the above system.

\begin{Remark}
{\rm
Under the assumptions of point 1 of Remark \ref{rkwave1},
the   observability estimate (\ref{OEelasticityint}) is obtained in the proof of Theorem 3.1 of \cite{Calva_98} (estimate (3.2) of \cite{Calva_98}).
\hspace{2 cm}\qed}
\end{Remark}

\subsection{Internal and boundary dampings}

Under the assumptions of subsection \ref{Wavebound} we consider the following elastodynamic system

 \begin{eqnarray}
\nonumber
& &u_{tt}(x,t) -\mu \Delta  u (x,t)-(\lambda+\mu) \nabla \hbox{ div }u
+b(x) 
u_t(x,t-\tau)=0\quad \mbox{\rm in}\quad\Omega\times
(0,+\infty),\label{P.1elaB}\\
& &u(x,t)=0, \quad x\in  \Gamma_1,\ t> 0\label{P.2elaB}\\
& &\sigma(u(x,t))\cdot \nu(x)=- k(x) u_t(x,t), \quad x\in  \Gamma_0,\ t> 0\label{P.3elaB}\\
& &u(x,0)=u_0(x)\quad \mbox{\rm and}\quad u_t(x,0)=u_1(x)\quad \hbox{\rm
in}\quad\Omega ,\label{P.4elaB}\\
& &\sqrt{b} u_t(x,t)=f^0(t)\quad \mbox{\rm in} \ \omega_2\times (-\tau , 0), \label{P.5elaB}
\end{eqnarray}
with
initial
data $(u_0, u_1, f^0)\in H^1_0(\Omega)^n\times L^2(\Omega)^n\times L^2((-\tau ,0);L^2(\omega_2)^n)$
and
$$
\sigma(u)=\mu (\sum_{i=1}^n\partial_i(u_j) \nu_i)_{j=1}^n+(\lambda+\mu) (\hbox{ div }u)\, \nu \hbox{ on } \Gamma_0.
$$

This problem enters into our abstract setting, once we take
$
H=L^2(\Omega)^n$,
   $A$ defined in the previous subsection, and
$B_1$ and $B_2$ defined as in subsection  \ref{Wavebound}.

As  $B_1$ is not bounded, we need to assume that
there exists a time $\overline T>0$ such that for every time $T>\overline T$
there is a
constant $c,$ depending on $T$ but independent of the initial data, such that
\begin{equation}\label{OEelasticitybdy}
\frac12((w_0,w_0)_V+\int_\Omega |w_1|^2\,dx)\le c\int_0^{T}\int_{\Gamma_0}k(x) |w_t|^2(x,s) dx ds,
\end{equation}
for every weak solution $w$ of  the non delayed  system
\begin{eqnarray}
\nonumber
& &w_{tt}(x,t) -\mu \Delta  w (x,t)-(\lambda+\mu) \nabla \hbox{ div }w
=0\quad \mbox{\rm in}\quad\Omega\times
(0,+\infty),\label{P.1elaC}\\
& &w(x,t)=0, \quad x\in  \Gamma_1,\ t> 0\label{P.2elaC}\\
& &\sigma(w(x,t))\cdot \nu(x)=- k(x) w_t(x,t), \quad x\in  \Gamma_0,\ t> 0\label{P.3elaC}\\
& &w(x,0)=w_0(x)\quad \mbox{\rm and}\quad w_t(x,0)=w_1(x)\quad \hbox{\rm
in}\quad\Omega ,\label{P.4elaC},
\end{eqnarray}
for initial
data $(w_0, w_1)\in H^1_0(\Omega)^n\times L^2(\Omega)^n$.

Again if such an estimate holds, the stability result from section \ref{st} can be applied to the   system
(\ref{P.1elaB})--(\ref{P.5elaB}).

\begin{Remark}
{\rm
1. Under the assumptions of point 1 of Remark \ref{rkwavebdy},
the   observability estimate (\ref{OEelasticitybdy}) is   proved in  \cite{Bey:03}.
\\
2. If we assume that the boundary of $\Omega$ is smooth and that
$$
(x-x_0)\cdot \nu(x)\leq 0 \hbox{ on } \Gamma_1,
$$
then the  observability estimate (\ref{OEelasticitybdy}) is   proved in Lemma 3.2 of \cite{horn}.
}\hspace{2 cm}\qed
\end{Remark}

\section{The Petrovsky system\label{PP}}

\hspace{5mm}

\setcounter{equation}{0}

\subsection{Hinged boundary conditions}

Let $\Omega\subset\RR^n$ be an open bounded set   with a boundary
$\partial\Omega$ of class $C^4$ (as before this regularity could be weakened).

Let us consider the initial boundary value  problem

 \begin{eqnarray}
& &u_{tt}(x,t) +\Delta^2 u (x,t)+b_1(x)  u_t(x,t)+b_2(x) 
u_t(x,t-\tau)=0\ \mbox{\rm in}\ \Omega\times
(0,+\infty),\label{P..1}\\
& &u (x,t) =\Delta u (x,t)=0\quad \mbox{\rm on}\quad\partial\Omega\times
(0,+\infty),\label{P..2}\\
& &u(x,0)=u_0(x)\quad \mbox{\rm and}\quad u_t(x,0)=u_1(x)\quad \hbox{\rm
in}\quad\Omega,\label{P..3}\\
& & u_t(x,t)=f^0(x,t) \quad \mbox{\rm in}\ \omega_2\times (-\tau ,0),\label{IP}
\end{eqnarray}
with
initial
data $(u_0, u_1, f^0)\in (H^2(\Omega)\cap H^1_0(\Omega))\times L^2(\Omega)\times L^2((-\tau ,0);L^2(\omega_2))$
and
$b_1,b_2$ satisfying
 the same assumptions as in subsection \ref{WWint}.

Now, we take
$H=L^2(\Omega)$ and let $A$ be the operator
\be\label{defbihar}
A:{\mathcal D}(A)\rightarrow H: \quad\quad u\rightarrow \Delta^2 u,
\ee
where
$${\mathcal D}(A)=\{v\in H^1_0(\Omega)\cap
H^4(\Omega): \Delta u=0 \hbox{ on } \partial\Omega\}.
$$

The operator $A$ is   self--adjoint and positive, has a compact inverse in $H$
and satisfies ${\mathcal D}(A^{1/2})=H^2(\Omega)\cap H^1_0(\Omega).$
We then define $U_i=L^2(\omega_i)$ and the operators $B_i, i=1,2,$ by (\ref{defBi}).
So, problem (\ref{P..1})--(\ref{IP}) enters in the abstract framework
(\ref{1.1})--(\ref{1.3}).

If an observability estimate of the associated conservative system holds,
then the results of section \ref{st} apply also to the plate model.

\begin{Remark}
{\rm
 Under the assumptions of point 1 of Remark \ref{rkwave1} and the additional regularity of the boundary,
it is well--known that an observability estimate of the associated conservative system holds, see Proposition 7.5.7 (see also Example 11.2.4) of \cite{TucsnakWeiss}.
}\hspace{2 cm}\qed
\end{Remark}

\subsection{Clamped boundary conditions}

Let $\Omega\subset\RR^2$ be an open bounded set   with a boundary
$\partial\Omega$ of class $C^4$.

Here we   consider the initial boundary value  problem

 \begin{eqnarray}
& &u_{tt}(x,t) +\Delta^2 u (x,t)+b_1(x)  u_t(x,t)+b_2(x) 
u_t(x,t-\tau)=0\ \mbox{\rm in}\ \Omega\times
(0,+\infty),\label{P.1Dir}\\
& &u (x,t) =\frac{\partial u}{\partial\nu} (x,t)=0\quad \mbox{\rm on}\quad\partial\Omega\times
(0,+\infty),\label{P.2Dir}\\
& &u(x,0)=u_0(x)\quad \mbox{\rm and}\quad u_t(x,0)=u_1(x)\quad \hbox{\rm
in}\quad\Omega,\label{P.3Dir}\\
& & u_t(x,t)=f^0(x,t) \quad \mbox{\rm in}\ \omega_2\times (-\tau ,0),\label{IPP}
\end{eqnarray}
with
initial
data $(u_0, u_1, f^0)\in H^2_0(\Omega)\times L^2(\Omega)\times L^2((-\tau ,0);L^2(\omega_2))$ where
$$H^2_0(\Omega ):=\{\ \varphi\in H^2(\Omega)\,:\  u =\frac{\partial u}{\partial\nu}=0\ \mbox{\rm on}\quad\partial\Omega\,\},$$
and
$b_1,b_2$ satisfy the same assumptions than in the previous subsection.

As usual, if an observability estimate holds for the associated conservative system,
 the results of section \ref{st} can be applied to this model.

\begin{Remark}
{\rm
 Under the assumptions of point 1 of Remark \ref{rkwave1} and the additional assumptions of this subsection, the observability estimate for the associated conservative system has been recently proved by the authors (see \cite{NPADE12}, Theorem 6.1).
}\hspace{2 cm}\qed
\end{Remark}

\bigskip

 {\em E-mail address,}
\\
\quad Serge Nicaise: \quad{\tt \bf
snicaise@univ-valenciennes.fr}
\\
\quad Cristina Pignotti: \quad{\tt \bf
pignotti@univaq.it}

\end{document}